\documentclass[12pt]{article}

\begin{document}
\thispagestyle{empty}
\newcommand{\tk}{2^{\textstyle \kappa}}
\newcommand{\tkk}{2^{\textstyle{2^{\textstyle \kappa}}}}
\newcommand{\ka}{\kappa^{\textstyle \ast}}
\newcommand{\kaa}{\kappa^{\textstyle \ast\ast}}
\newcommand{\kd}{\kappa^{\textstyle \diamond}}
\newcommand{\kdd}{\kappa^{\textstyle \diamond\diamond}}
\newcommand{\So}{S_{\textstyle 1}}
\newcommand{\St}{S_{\textstyle 2}}
\newcommand{\Bo}{B_{\textstyle 1}}
\newcommand{\Bt}{B_{\textstyle 2}}
\newcommand{\RB}{R_{\textstyle B}}
\newcommand{\fB}{f_{\textstyle B}}
\newcommand{\idk}{{\rm id}_{\textstyle \kappa}}
\newcommand{\RP}{R_{\textstyle \Psi}}
\newcommand{\Cn}{{\cal C}_{\textstyle n}}
\newcommand{\R}{{\cal R}}
\def\thefootnote{}
\centerline{\large \bf A stronger form of the theorem constructing}
\centerline{\large\bf a rigid binary relation on any set}
\vskip 2.0truecm
\centerline{\large \bf Apoloniusz Tyszka}
\vskip 2.0truecm
\par
{\bf Summary.} On every set $A$ there is a rigid binary relation
i.e. such a relation $\R \subseteq A \times A$ that there is no
homomorphism $\langle A,\R \rangle \rightarrow \langle A,\R \rangle$
except the identity (Vop{\v{e}}nka et al. [1965]). We prove that for
each infinite cardinal number $\kappa$ if card~$A \leq \tk$, then
there exists a relation $\R \subseteq A \times A$ with the following
property:
\\
\\
\centerline{
$\forall x \in A$
$\exists^{\textstyle
{\{x\} \subseteq A(x) \subseteq A}}_{\textstyle
{{\rm card}~A(x) \leq \kappa}}$
$\forall^{\textstyle
{f:A(x) \rightarrow A}}_{\textstyle
{f \neq {\rm id}_{\textstyle
A(x)}}}$ $f$ {\rm is not a homomorphism of} $\R$}
\vskip 0.4truecm
\par
\noindent
which implies that $\R$ is rigid. If a relation
$\R \subseteq A \times A$ has the above property,
then ${\rm card}~A \leq \tk$.
\vskip 2.0truecm
\par
On every set $A$ there is a rigid binary relation, i.e. such a
relation $\R \subseteq A \times A$ that there is no homomorphism
$\langle A,\R \rangle \rightarrow \langle A,\R \rangle$ except
the identity
(\cite{nesetril1},\cite{nesetril2} \cite{pultr},\cite{vopenka}).
Conjectures~1 and 2 below strengthen this theorem.
\newpage
{\bf Conjecture~1} (\cite{tyszka1},\cite{tyszka2}).
If $\kappa$ is an infinite cardinal number and
${\rm card}~A \leq \tkk$, then there exists a relation
$\R\subseteq A\times A$ which satisfies the following
condition~$(\ka)$:
\vskip 0.3truecm
\begin{description}
\item{$(\ka)$}
\centerline{$\forall^{\textstyle x,y \in A}_{\textstyle x \neq y}$
$\exists^{\textstyle
{\{x\} \subseteq A(x,y) \subseteq A}}_{\textstyle
{{\rm card}~A(x,y) \leq \kappa}}$
$\forall^{\textstyle
{f:A(x,y) \rightarrow A}}_{\textstyle
{f(x)=y}}$}
\\
$f$ {\rm is not a homomorphism of} $\R$.
\end{description}
\par
{\bf Proposition~1a} (\cite{tyszka2}).
If $\kappa$ is an infinite cardinal number,
$\R \subseteq A \times A$ satisfies condition~$(\ka)$ and
${\rm card}~\widetilde{A} \leq {\rm card}~A$, then there exists a
relation
$\widetilde{\R} \subseteq \widetilde{A} \times \widetilde{A}$
which satisfies condition~$(\ka)$.
\vskip 0.2truecm
\par
{\bf Proposition~1b} (\cite{tyszka1}). If $\R \subseteq A \times A$
satisfies condition~$(\ka)$, then $\R$ is rigid.
If $\kappa$ is an infinite cardinal number and a relation
$\R\subseteq A\times A$ satisfies condition~$(\ka)$,
then ${\rm card}~A \leq \tkk$.
\vskip 0.2truecm
\par
{\bf Theorem~1} (\cite{tyszka2}). Conjecture 1 is valid for
$\kappa=\omega$.
\vskip 0.2truecm
\par
{\bf Conjecture~2} (\cite{tyszka1},\cite{tyszka2}).
If $\kappa \neq 0$ is a limit cardinal number and
${\rm card}~A \leq 2^{\textstyle
{\rm sup} \{2^{\textstyle \alpha}:
\alpha \in {\rm Card}, \alpha < \kappa \}}$,
then there exists a relation $\R \subseteq A\times A$ which
satisfies the following condition~$(\kaa)$:
\begin{description}
\item{$(\kaa)$}
\centerline{
$\forall^{\textstyle x,y \in A}_{\textstyle x \neq y}$
$\exists^{\textstyle
{\{x\} \subseteq A(x,y) \subseteq A}}_{\textstyle
{{\rm card}~A(x,y) < \kappa}}$
$\forall^{\textstyle
{f:A(x,y) \rightarrow A}}_{\textstyle
{f(x)=y}}$}
\\
{$f$ {\rm is not a homomorphism of} $\R$.}
\end{description}
\par
{\bf Proposition~2a} (\cite{tyszka2}).
If $\kappa \neq 0$ is a limit cardinal number,
$\R \subseteq A \times A$ satisfies condition~$(\kaa)$
and ${\rm card}~\widetilde{A} \leq {\rm card}~A$,
then there exists a relation
$\widetilde{\R} \subseteq \widetilde{A} \times \widetilde{A}$
which satisfies condition~$(\kaa)$.
\vskip 0.2truecm
\par
{\bf Proposition~2b} (\cite{tyszka1}). If $\R \subseteq A \times A$
satisfies condition~$(\kaa)$, then $\R$ is rigid.
If $\kappa \neq 0$ is a limit cardinal number and a relation
$\R\subseteq A\times A$ satisfies condition~$(\kaa)$,
then ${\rm card}~A \leq 2^{\textstyle
{\rm sup} \{2^{\textstyle \alpha}:
\alpha \in {\rm Card}, \alpha < \kappa \}}$.
\vskip 0.2truecm
\par
{\bf Theorem~2} (\cite{tyszka1},\cite{tyszka2}).
Conjecture 2 is valid for $\kappa=\omega$.
\vskip 0.2truecm
\par
In this article we prove a changed form of Conjecture 1 which holds
for all infinite cardinal numbers $\kappa$, see
Theorems~3 and 4.
\vskip 0.2truecm
\par
{\bf Theorem~3.} If $\kappa$ is an infinite cardinal number and
${\rm card}~A \leq \tk$, then there exists a relation
$\R \subseteq A \times A$ which satisfies the following
condition~$(\kd)$:
\vskip 0.3truecm
\begin{description}
\item{$(\kd)$}
\centerline{
$\forall x \in A$
$\exists^{\textstyle
{\{x\} \subseteq A(x) \subseteq A}}_{\textstyle
{{\rm card}~A(x) \leq \kappa}}$
$\forall^{\textstyle
{f:A(x) \rightarrow A}}_{\textstyle
{f \neq {\rm id}_{\textstyle
A(x)}}}$}
\\
$f$ {\rm is not a homomorphism of} $\R$.
\end{description}
\par
{\it Proof.} It is known
(\cite{hedrlin},\cite{nesetril1},\cite{pultr})
that for each infinite cardinal number $\kappa$ there exists a rigid
symmetric relation $R \subseteq \kappa \times \kappa$. Let $\Phi$
denote the family of all relations
$S \subseteq \kappa \times \kappa$ which satisfy:
\par
(1) $R \subseteq S$,
\par
(2) for each $\alpha, \beta \in \kappa$ if $\alpha \neq \beta$,
then $\alpha S \beta$ or $ \beta S \alpha$,
\par
(3) for each $\alpha, \beta \in \kappa$ if $\alpha S \beta$
and $\beta S \alpha$, then $ \alpha R \beta$ and
$\beta R \alpha$.
\vskip 0.2truecm
\par
Since $R$ is rigid
\par
(4) $R \subseteq \{(\alpha,\beta):
\alpha,\beta \in \kappa, \alpha \neq \beta\}$.
\vskip 0.2truecm
\par
By (1) and (3) the following Lemma 1 holds true.
\vskip 0.2truecm
\par
{\bf Lemma~1.} If $\So,\St \in \Phi$ and
$f:\langle\kappa,\So\rangle \rightarrow \langle\kappa,\St\rangle$
is a homomorphism, then
$f:\langle\kappa,R\rangle \rightarrow \langle\kappa,R \rangle$ is
a homomorphism.
\vskip 0.2truecm
\par
{\bf Lemma~2.} For every $\So, \St \in \Phi$ if $\So \neq \St$, then
$\idk: \langle \kappa,
\So \rangle \rightarrow \langle \kappa,\St \rangle$
is not a homomorphism.
\par
{\it Proof.} Applying (3) and (4) we obtain two cases.
First case:
there exist $\alpha,\beta \in \kappa$, $\alpha \neq \beta$ such
that $(\alpha,\beta) \in \So$ and $(\alpha,\beta) \not\in \St$,
so $\idk$ is not a homomorphism.
Second case:
there exist $\alpha,\beta \in \kappa$, $\alpha \neq \beta$
such that $(\alpha,\beta) \in \St$ and $(\alpha,\beta) \not\in \So$.
By (2) $(\beta,\alpha) \in \So$. It suffices to prove that
$(\beta,\alpha) \not\in \St$. Suppose, on the contrary, that
$(\beta,\alpha) \in \St$. By (3) $(\alpha,\beta) \in R$,
so by (1) $(\alpha,\beta) \in \So$, a contradiction.
\vskip 0.2truecm
\par
{\bf Lemma~3.} ${\rm card}~\Phi= \tk$.
\par
{\it Proof.} Let $T:=\{\{\alpha,\beta\}:
\alpha,\beta~\in \kappa, \alpha \neq \beta,
(\alpha,\beta) \not\in R\}$.
It suffices to prove that ${\rm card}~T=\kappa$.
Suppose, on the contrary, that ${\rm card}~T < \kappa$.
Hence ${\rm card} \bigcup T < \kappa$
and consequently ${\rm card}~(\kappa \setminus \bigcup T) = \kappa$.
For each $\alpha,\beta \in \kappa \setminus \bigcup T$ if
$\alpha \neq \beta$, then $(\alpha,\beta) \in R$. From this and (4)
any non-identical injection from
$\kappa$ into $\kappa \setminus \bigcup T$ is a homomorphism of
$R$. This contradiction completes the proof of Lemma~3.
\vskip 0.2truecm
\par
Now we turn to the main part of the proof. For each
$\emptyset \neq \Psi \subseteq \Phi$ we define the relation
$\RP \subseteq (\kappa \times \Psi) \times (\kappa \times \Psi)$
by the following formula:
\\
\centerline{$\forall \alpha,\beta \in \kappa
\forall \So, \St \in \Psi
\Bigl(((\alpha,\So),(\beta,\St) ) \in \RP
\Longleftrightarrow (\alpha,\beta) \in \So=\St\Bigr)$.}
\par
\noindent
In other words, the graph corresponding to the relation $\RP$
is a disjoint union of graphs belonging to $\Psi$.
By Lemma~3 it suffices to prove that $\RP$ satisfies
condition~$(\kd)$. Let $(\lambda,\So) \in \kappa \times \Psi$.
We prove that
$(\kappa \times \Psi)((\lambda,\So)):=
\kappa~\times~\{\So\}$ satisfies condition~$(\kd)$.
\par
\vskip 0.2truecm
Suppose, on the contrary, that
$f:\kappa\times\{\So\} \rightarrow \kappa \times \Psi$
is a homomorphism of
$\RP$ and $f \neq {\rm id}_{\textstyle \kappa \times \{\So\}}$.
Then there exist $\alpha, \beta \in \kappa$ and $\St \in \Psi$
such that $f((\alpha,\So))=(\beta,\St)$ and
$(\alpha,\So) \neq (\beta,\St)$.
By (2) for each $\gamma \in \kappa \setminus \{\alpha\}$
$\alpha \So \gamma$ or $\gamma \So \alpha$. From this for each
$\gamma \in \kappa \setminus \{\alpha\}$
$(\alpha, \So) \RP (\gamma,\So)$ or $(\gamma,\So) \RP (\alpha,\So)$.
Therefore $f((\alpha,\So)) \RP f((\gamma,\So))$ or \hspace{0.1truecm}
$f((\gamma,\So)) \RP f((\alpha,\So))$ \hspace{0.1truecm} and
\hspace{0.1truecm}
consequently \hspace{0.1truecm} $(\beta,\St) \RP f((\gamma,\So))$
\hspace{0.1truecm} or $f((\gamma,\So)) \RP (\beta,\St)$.
In both cases there exists a
$\delta \in \kappa$ such that $f((\gamma,S_{1}))=(\delta,\St)$.
It implies that $f$ maps $\kappa \times \{\So\}$ into
$\kappa \times \{\St\}$. Let $\pi: \{\So\} \rightarrow \{\St\}$.
There is a uniquely determined transformation
$\widetilde{f}: \kappa \rightarrow \kappa$ such that
$f$ $=$ $\langle\widetilde{f},\pi\rangle$.
Obviously, $\widetilde{f}(\alpha)=\beta$ and
$\widetilde{f}:
\langle\kappa,\So\rangle \rightarrow \langle \kappa,\St\rangle$
is a homomorphism. By Lemma~1
$\widetilde{f}:
\langle\kappa, R\rangle \rightarrow \langle\kappa, R\rangle$
is a homomorphism. Since $R$ is rigid
$\widetilde{f}=\idk$. Therefore
$\alpha=\widetilde{f}(\alpha)=\beta$ and
$\idk:\langle\kappa,\So\rangle \rightarrow \langle\kappa,\St\rangle$
is a homomorphism. On the other hand, $\alpha=\beta$
and $(\alpha,\So) \neq (\beta,\St)$ implies $\So \neq \St$.
It is impossible by Lemma~2. This contradiction completes the
proof of Theorem~3.
\vskip 0.2truecm
\par
{\bf Remark.} It is easy to observe that
condition~$(\kd)$ implies condition~$(\ka)$.
Obviously, if $\R \subseteq A \times A$ satisfies condition~$(\kd)$,
then $\R$ is rigid.
\vskip 0.2truecm
\par
We will show an alternative, algebraic method for proving Theorem~3.
This method described in the proof of Theorem~4 gives more, namely
a symmetric relation satisfying condition~$(\kd)$.
Unfortunately, in contradistinction to the relation constructed in
the proof of Theorem~3, a direct description of such a relation is
very complicated.
\vskip 0.2truecm
\par
{\bf Theorem~4.} If $\kappa$ is an infinite cardinal number and
${\rm card}~A \le \tk$, then there exists a symmetric relation
$\R \subseteq A \times A$ which satisfies condition~$(\kd)$.
\par
{\it Proof.} Following \cite{pultr} let {\bf Graph} denote the
category of graphs and their homomorphisms. The objects of
{\bf Graph} are couples $(X,R)$ with $R \subseteq X \times X$,
the morphisms from $(X,R)$ to $(X',R')$ are triples
$((X',R'),f,(X,R))$ with $f: X \rightarrow X'$ such that
$(f(x),f(y)) \in R'$ whenever $(x,y) \in R$, and it is viewed
as a concrete category endowed with the natural forgetful functor.
Let $\Cn$ $(n \ge 3)$ denote the category of connected $n$-chromatic
undirected graphs and their homomorphisms. It is known
(see \cite{pultr}, Theorem 4.12 on page 113), that for every
$n \ge 3$ there is a strong embedding
$F: {\bf Graph} \rightarrow \Cn$ that transforms objects of
the cardinality $\kappa$ into objects of the cardinality $\kappa$.
\par
Let $R \subseteq \kappa \times \kappa$ be a rigid symmetric
relation (undirected graph). Considering all possible orientations
of $R$ we obtain $\tk$ rigid graphs with the property that there
is no homomorphism between any two distinct graphs. Using the strong
embedding $F: {\bf Graph} \rightarrow \Cn$ $(n \ge 3)$
constructed in \cite{pultr} we obtain a family $\Gamma$ of $\tk$
rigid undirected graphs on $\kappa$ with the property that there is
no homomorphism between any two distinct graphs. Since $\Gamma$
consists of connected graphs, for every $\Delta \subseteq \Gamma$
a disjoint union of graphs belonging to $\Delta$ satisfies
condition~$(\kd)$, the proof is similar to the proof that $\RP$
satisfies condition~$(\kd)$. This completes the proof of Theorem 4.
\newpage
{\bf Theorem~5.} If $\kappa$ is an infinite cardinal number and
a relation $\R \subseteq A \times A$ satisfies condition~$(\kd)$,
then ${\rm card}~A \leq \tk$.
\par
{\it Proof.} Suppose, on the contrary, that
$\R \subseteq A \times A$ satisfies condition~$(\kd)$ and
${\rm card}~A > \tk$.
For each $x \in A$ we choose the set $A(x)$ from condition~$(\kd)$
in such a way that ${\rm card}~A(x)=\kappa$. Let
${\cal B}:=\{A(x): x \in A\}$. Since $\bigcup{\cal B}=A$ we conclude
that ${\rm card}~{\cal B}={\rm card}~A$. For each $B \in {\cal B}$
we choose a bijective $\fB: \kappa \rightarrow B$ and define the
relation $\RB \subseteq \kappa \times \kappa$ by the following
formula:
\\
\centerline{
$\forall \alpha,\beta \in \kappa$
$\Bigl((\alpha,\beta) \in \RB \Longleftrightarrow
(\fB(\alpha),\fB(\beta)) \in \R\Bigr)$.}
\par
\noindent
Let ${\cal B} \ni B \stackrel{{\textstyle h}}{\longrightarrow}
\RB \in {\cal P}(\kappa \times \kappa)$.
Since ${\rm card}~{\cal B}=
{\rm card}~A > \tk = {\rm card}~{\cal P}(\kappa \times \kappa)$
we conclude that there exist $\Bo,\Bt \in {\cal B}$ such that
$\Bo \neq \Bt$ and $h(\Bo)=h(\Bt)$. Hence
$\langle \Bo,\R \rangle \stackrel{{\textstyle f_{\textstyle \Bt}
\circ (}{\textstyle f_{\textstyle {\Bo}})^{{\textstyle -1}}}}
{\longrightarrow} \langle \Bt,\R \rangle$ is a non-identical
isomorphism. This contradiction completes the proof.
\vskip 0.2truecm
\par
{\bf Conjecture~3.} If $\kappa \neq 0$ is a limit cardinal number
and ${\rm card}~A \leq {\rm sup} \{2^{\textstyle \alpha}:
\alpha \in {\rm Card}, \alpha < \kappa \}$,
then there exists a relation
$\R \subseteq A \times A$ which satisfies the following
condition~$(\kdd)$:
\vskip 0.3truecm
\begin{description}
\item{$(\kdd)$}
\centerline{
$\forall x \in A$
$\exists^{\textstyle
{\{x\} \subseteq A(x) \subseteq A}}_{\textstyle
{{\rm card}~A(x) < \kappa}}$
$\forall^{\textstyle
{f:A(x) \rightarrow A}}_{\textstyle
{f \neq {\rm id}_{\textstyle
A(x)}}}$}
\\
$f$ {\rm is not a homomorphism of} $\R$.
\end{description}
\par
{\bf Theorem~6.} If $\kappa \neq 0$ is a limit cardinal number and
a relation $\R \subseteq A \times A$ satisfies condition~$(\kdd)$,
then ${\rm card}~A \leq {\rm sup} \{2^{\textstyle \alpha}:
\alpha \in {\rm Card}, \alpha < \kappa \}$.
\par
{\it Proof.} Suppose, on the contrary, that
$\R \subseteq A \times A$ satisfies condition~$(\kdd)$ and
${\rm card}~A > {\rm sup} \{2^{\textstyle \alpha}:
\alpha \in {\rm Card}, \alpha < \kappa \}$.
For each $x \in A$ we choose the set $A(x)$ from condition~$(\kdd)$.
Let ${\cal B}:=\{A(x): x \in A\}$. Since $\bigcup{\cal B}=A$ we conclude
that ${\rm card}~{\cal B}={\rm card}~A$. For each $B \in {\cal B}$
we choose a bijective $\fB: {\rm card}~B \rightarrow B$ and define
the relation $\RB \subseteq {\rm card}~B \times {\rm card}~B$ by the
following formula:
\\
\centerline{
$\forall \alpha,\beta \in {\rm card}~B$
$\Bigl((\alpha,\beta) \in \RB \Longleftrightarrow
(\fB(\alpha),\fB(\beta)) \in \R\Bigr)$.}
\par
\noindent
Let
\begin{displaymath}
{\cal B} \ni B \stackrel{{\textstyle h}}{\longrightarrow}
({\rm card}~B, \RB) \in \bigcup_{\stackrel
{\textstyle \alpha \in {\rm Card}}
{\textstyle \alpha < \kappa}}
\{\alpha\} \times {\cal P}(\alpha \times \alpha).
\end{displaymath}
Since
\begin{displaymath}
{\rm card}~{\cal B}=
{\rm card}~A > {\rm sup} \{2^{\textstyle \alpha}:
\alpha \in {\rm Card}, \alpha < \kappa \}=
{\rm card}~\bigcup_{\stackrel
{\textstyle \alpha \in {\rm Card}}
{\textstyle \alpha < \kappa}}
\{\alpha\} \times {\cal P}(\alpha \times \alpha)
\end{displaymath}
we conclude that there exist $\Bo,\Bt \in {\cal B}$ such that
$\Bo \neq \Bt$ and $h(\Bo)=h(\Bt)$. Hence
$\langle \Bo,\R \rangle \stackrel{{\textstyle f_{\textstyle \Bt}
\circ (}{\textstyle f_{\textstyle {\Bo}})^{{\textstyle -1}}}}
{\longrightarrow} \langle \Bt,\R \rangle$ is a non-identical
isomorphism. This contradiction completes the proof.
\vskip 0.2truecm
\par
{\bf Proposition~3.} Obviously, if $\R \subseteq A \times A$
satisfies condition~$(\kdd)$, then $\R$ is rigid. By Theorems 3 and 6,
if $\kappa \neq 0$ is a limit cardinal number,
$\R \subseteq A \times A$ satisfies condition~$(\kdd)$
and ${\rm card}~\widetilde{A} \leq {\rm card}~A$,
then there exists a relation
$\widetilde{\R} \subseteq \widetilde{A} \times \widetilde{A}$
which satisfies condition~$(\kdd)$.
\vskip 0.2truecm
\par
{\bf Theorem 7.} Conjecture 3 is valid for $\kappa=\omega$ i.e.
there exists a relation $\R \subseteq \omega \times \omega$
satisfying condition~$(\omega^{\textstyle \diamond\diamond})$.
\par
{\it Proof.} The relation
$\R:=\{(i,i+1):
i \in \omega \} \cup \{(0,2)\} \subseteq \omega \times \omega$
satisfies condition~$(\omega^{\textstyle \diamond\diamond})$. Indeed,
for each $i \in \omega$ the set $A(i):=\{j \in \omega: j \leq i+2 \}$
is adequate for property~$(\omega^{\textstyle \diamond\diamond})$.
\vskip 0.2truecm
\par

\begin{flushleft}
Apoloniusz Tyszka\\
Technical Faculty\\
Hugo Ko{\l}{\l}\c{a}taj University\\
Balicka 104, 30-149 Krak\'{o}w, Poland\\
E-mail: {\it rttyszka@cyf-kr.edu.pl}\\
\end{flushleft}
\end{document}